**Continuous Synthesis of Diazo Acetonitrile: From Experiments to Physical and Grey-Box Modeling**

Dr. Marco Baldan[1,]*, Dr. Sebastian Blauth[1], Dr. Dušan Bošković[2], Dr. Christian Leithäuser[1], Dr. Alexander Mendl[2], Ligia Radulescu[2], Maud Schwarzer[2], Heinrich Wegner[2], and Prof. Dr. Michael Bortz[1]

Diazo compounds are gathering interest for their potential in promoting greener synthesis routes. We investigate, at a lab-scale, the continuous synthesis of diazo acetonitrile (DAN) using a micro-structured flow reactor and a flow reaction calorimeter. Data concerning DAN formation in the former, and relative to reaction heat and gas flow rate in the latter, are collected. We present both a physical and a grey-box simulation model, both of which are calibrated to our measurements. Both models provide valuable insights into the DAN synthesis. The grey-box approach is useful to incorporate the complex chemical reaction pathways for DAN synthesis and decomposition that are currently hard to address with the physical model.

**Keywords:**

**Author affiliations**

[1]Fraunhofer ITWM, Fraunhofer-Platz 1, 67663 Kaiserslautern, Germany.

[2]Fraunhofer ICT, Joseph-von-Fraunhofer-Straße 7, 76327 Pfinztal, Germany.

Email corresponding author: marco.baldan@itwm.fraunhofer.de

ORCID iDs of the authors

0000-0002-5803-3150 (Marco Baldan)

0000-0001-9173-0866 (Sebastian Blauth)

0000-0001-8169-2907 (Michael Bortz)

0009-0001-6563-6611 (Dušan Bošković)

0000-0001-8936-9805 (Christian Leithäuser)

0000-0003-1735-4684 (Alexander Mendl)

0000-0002-9917-1381 (Maud Schwarzer)

## 1  Introduction

Diazo compounds are increasingly gaining in importance in the synthesis of specialty chemicals, not only for their versatile reactivity but also due to their potential in promoting greener synthesis routes. Green chemistry principles aim to minimize the use of hazardous substances, reduce waste generation, and conserve energy during chemical processes. This seems paradoxical at first since diazo compounds themselves are considerably hazardous substances. And yet, through appropriate process technology, they can contribute significantly to the development of novel and innovative specialty chemicals.

One key advantage is their ability to enable selective and efficient transformations. Diazo compounds can undergo a wide range of reactions, such as cyclopropanation, cycloadditions, and C-H functionalization, which are often difficult to achieve using traditional methods [1]. An important aspect of diazo compounds in the synthesis of specialty chemicals is their contribution to atom efficiency. Atom efficiency is a key principle in green chemistry that focuses on maximizing the utilization of atoms in a reaction, minimizing waste, and improving overall resource utilization. Diazo compounds offer the advantage of enabling more direct synthesis routes, often requiring fewer synthetic steps compared to alternative methods [2]. This direct functionalization approach avoids the need for protecting groups or multiple reaction sequences that are often required in conventional synthesis methods. By eliminating unnecessary steps, the overall atom economy of the process is improved, resulting in a more sustainable and efficient synthesis route.

However, the handling and storage of diazo compounds can pose significant safety risks due to their instability and often explosive nature [1]. When it comes to the synthesis of specialty chemicals involving diazo compounds, continuous process technology provides unique benefits that enhance both process safety and efficiency. Instead of storing and handling pre-formed diazo compounds, these compounds are synthesized in-situ and when needed [3–6], which eliminates the need for storing potentially dangerous materials. With this in mind, numerous process examples for continuous diazo reagent synthesis have been developed in recent years, most notably of particularly useful and equally hazardous diazomethane [7–11] but also other compounds including ethyl diazoacetate [12, 13], difluoromethyl, and trifluoromethyl diazomethane [14–16].

Another compound that has received little attention so far but is likewise useful is diazo acetonitrile (DAN). It can be used to introduce nitrile groups, which are found in a great number of pharmaceutical and agrochemical active ingredients, into molecules. Nitrile compounds can also be reduced to amines, which also offers a wide range of applications [17]. DAN was synthesized for the first time in 1898, but presumably because of the potential dangers due to its high nitrogen content, which led to some serious accidents [18], it was not being utilized much. In recent years, Mykhailiuk [19] presented a procedure where DAN was formed in situ in a batch approach and directly reacted in the same flask with an alkyne to form different pyrazoles. Empel [20] took it a step further by forming DAN continuously in a micromixer and directly feeding this mixture in a standard reaction flask where subsequent insertion reactions were realized. In both cases, DAN was formed in aqueous phase and the subsequent reaction took place in a two-phase mixture. DAN was not isolated or rather transferred to an organic phase for subsequent utilization as it is frequently done in the case of other diazo reagents [3], e.g., ethyl diazoacetate [13]. To enable scalability and allow for versatile utilization of DAN, it would be beneficial to continuously prepare an organic solution of DAN, which can be flexibly utilized in a likewise continuous subsequent reaction step. To realize such a process, it is necessary to study the formation of DAN in detail. So far, almost no reaction or kinetic data are available.

A particularly useful technique for analyzing reaction conditions and kinetics is reaction calorimetry as it provides fast and direct measurement of heat released during a chemical reaction. Therefore, reaction calorimetry can provide insights into the kinetics, rate, and thermal stability of the reaction. Additionally, reaction calorimetry can help in optimizing reaction conditions, identifying potential hazards, and designing safe and efficient processes. In recent years, a new generation of reaction calorimeters based on micro structured flow reactors and microreactors has emerged [21–24]. The flow reactors used therein provide excellent heat transfer capabilities due to their high specific surface area, which enables precise temperature control. Typically, these systems achieve a spatially resolved measurement by placing multiple heat flow sensors along a reaction channel. One of the first flow reaction calorimeters was developed by Antes et al. [25, 26] and a further developed version is used in this work for analyzing DAN generation. It measures the heat release rate simultaneously at different positions along the reactor which on one hand allows to determine optimal reaction parameters such as temperature, residence time, and reactant concentrations, and on the other hand is used to develop and validate reaction rate models.

Mathematical modelling of reactors [27] could provide valuable insights into the process and be useful in other contexts, like model-based optimization [28], model-based design of experiments [29] or model predictive control [30], to name a few. However, to the best of our knowledge, little literature is available regarding the modelling of systems containing diazo compounds, including DAN. Motivated by this, on the one hand, we consider a physical model for the chemically reactive flow [27]. On the other hand, we employ a grey-box plug-flow reactor (PFR) model of the reactor that incorporates a data-driven component, namely the expression for the source term of a reaction. Grey-box modelling represents an established alternative to purely physical approaches, particularly useful when the modelled phenomena either are not fully understood or first principle models are not accurate [31]. Grey-box modelling has been already applied to chemical reactors [32]. However, compared to physical models, grey-box models introduce a higher number of parameters, making them more prone to overfitting in the parameter identification [33].

This paper is structured as follows: after a brief description of the chemical reactions involved in the systems, Sections 2.2 and 2.3 describe the two lab-scale setups used for the continuous synthesis of DAN. The physical and the grey-box PFR model are then explained in Section 2.4 and 2.5, respectively. The kinetic identification problem is presented in Section 2.6. Finally, Section 3 investigates the results of the experiments together with our parameter identification approach. Section 4 concludes the paper.

## 2 Methods and Models

We now describe the experimental setups, the so-called "mixer" and "calorimeter" setups. Further, we introduce two mathematical models, namely a physical and a grey-box PFR model. Additionally, we present the parameter identification problems used to fit our models to the experimental data.

### 2.1 Reaction Model for Diazo Acetonitrile Synthesis and Decomposition

The diazo acetonitrile (DAN) is obtained from aminoacetonitrile hydrochloride ($H_2NCH_2CN\cdot HCl$, CAS: 6011-14-9) and sodium nitrite ($NaNO_2$, CAS: 7632-00-0) through the diazotization reaction [19]

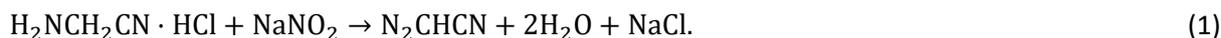

$$H_2NCH_2CN \cdot HCl + NaNO_2 \rightarrow N_2CHCN + 2H_2O + NaCl. \tag{1}$$

Additionally, due to its instability notably in aqueous solution at elevated temperature, also decomposition of DAN takes place forming nitrogen and a carbene

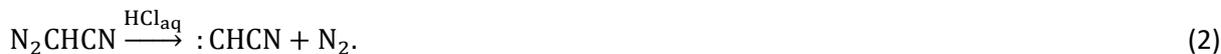

$$N_2CHCN \xrightarrow{HCl_{aq}} :CHCN + N_2. \tag{2}$$

These are the primary reactions taking place. The formed carbene can undergo numerous secondary reactions, forming dimers, oligomers, or an alcohol by the reaction with water molecules [34, 35]. These reactions, due to growing complexity, are not considered here.

### 2.2 Mixer: Setup for the Online FTIR Measurement

Similar to [6], a continuous flow setup is used to run the experiments. In our setup (Figure 1, Setup 1), the educts are fed by continuous operating syringe pumps (Syrdos, Hitec Zang GmbH), with aminoacetonitrile hydrochloride at a molar concentration of 2 mol/L and $NaNO_2$ at molar concentration of 2.4 mol/L. The pressure inside the reactor is controlled via a backpressure regulator (V-1, Zaiput Flow Technologies) and measured before entering a temperature bath (P-01).

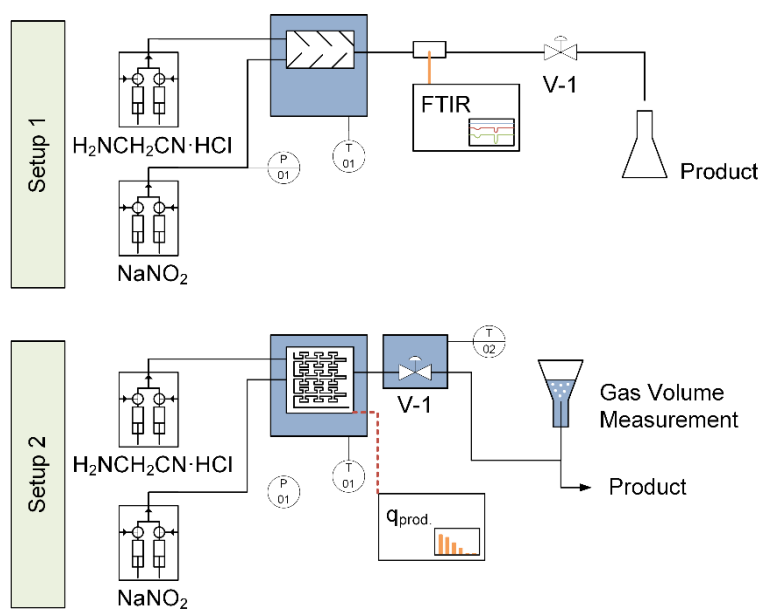

**Figure 1.** Reaction scheme of the two setups (setup 1: "mixer", setup 2: "calorimeter").

A glass reactor with a micro channel having a hydraulic diameter of 500 µm and a chaotic mixing structure with a volume of 1.4 ml (ST design, Little Things Factory) is used. The product mixture is analyzed by using a flow cell (S-PACT GmbH) equipped with an ATR-Probe and a FTIR analyzer (Matrix, Bruker Corporation).

Prior to the experiments, the system is flushed with deionized water and air bubbles were removed. The temperature is varied in a range from 20 °C to 70 °C and the flow rate in a range of 0.4 mL/min to 4 mL/min, keeping the nitrite access at 1.2. Experiments are conducted at a pressure of 1.4 bar gauge.

After reaching steady state conditions, the product is collected at the outlet and the online FTIR measurement is started. The product DAN shows an absorption band of the N=N stretching vibration with the maximum at 2110 cm$^{-1}$ that is not affected by other absorptions in the product mixture. As the mixture also contains gaseous decomposition products, signals are acquired at a high scan rate. The gas phase signals, indicated by a shifting baseline or by losing intensity in the diazo band, are then removed and the resulting spectra are averaged. The band of the diazo mode is integrated in the range from 2082 cm$^{-1}$ and 2135 cm$^{-1}$. Because DAN cannot be isolated due to its instability, no calibration can be done. Consequently, the band area is used as a relative value for the DAN content. The collected data are available in the supporting information.

### 2.3 Calorimeter: Setup for Measuring the Heat of Reaction and Gas Volume

The heat flow of reaction is measured in a continuously operating reaction calorimeter that has been previously developed at Fraunhofer ICT (Figure 1, Setup 2) [26]. Gas volumes are measured at ambient pressure at the outlet of the flow reactor to determine the degree of decomposition. To prevent the reactions from proceeding any further, the mixture is cooled in an ice bath (T-02 in Figure 1).

The flow reaction calorimeter contains a reaction plate with a meandering microchannel, enabling convective mixing due to alternating dean vortices [36]. The channel has a square shaped geometry with an edge length of 500 µm and a total volume of 889 µL. Micro-milling of stainless steel and a bonding process with a soldering foil is used to build the reaction plate at Fraunhofer IMM (Mainz, Germany). The reaction plate is divided into seven distinct zones (Figure 2). The heat flow in each zone is measured independently allowing space resolved heat flow measurements. This should ease the identification of the reaction kinetics and provide immediate detection of hot spots.

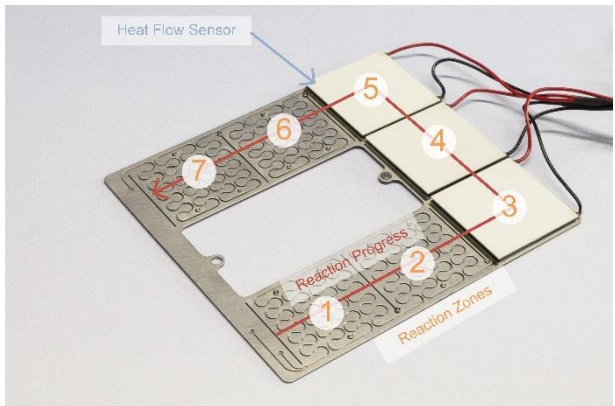

**Figure 2.**    Scheme of the reaction plate with the measuring zones (1 to 7, from inlet to outlet).

To measure the heat flow $q$, Seebeck elements are used on the top and the bottom of each reaction zone. Individual sensors are calibrated via an integrated heating foil [37]. The reaction plate equipped with the heat flow sensors is installed in a temperature bath to ensure isothermal conditions. To ensure that tempering the reactants does not influence the measurement, the reactants are passed through a pre-heating loop before entering the reaction plate (Figure 3). Technical efforts were made to minimize losses and convective induced heat. LabView® is used to read, process, and display the data, as well as to control the entire system including pumps, sensors, and thermostats.

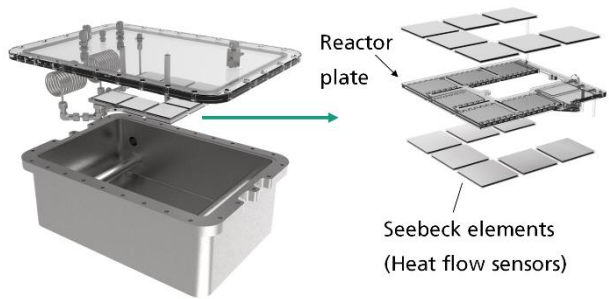

**Figure 3.**    Scheme of the temperature bath and the assembly of the heat flow sensors on the reaction plate.

The experiments are conducted at the same temperature span as in the mixer but at slightly different flow ranges (0.5 to 2 mL/min) and at several pressure levels (0 bar, 3 bar, and 6 bar, gauge pressure). After reaching steady state conditions, the heat flow of each zone is measured. At the same time, the gas volume measurement is started by measuring of the time for sampling of 5 mL of gaseous reaction products. Results are presented in the supporting information.

### 2.4    Physical Two-Phase Model for the Chemically Reactive Flow

We now describe the physical model of the DAN synthesis in a continuous flow setting. Let $\Omega$ be the domain of the reactor. Its boundary $\Gamma$ is divided into the inlet $\Gamma_{\text{in}}$, wall boundary $\Gamma_{\text{wall}}$, and outlet $\Gamma_{\text{out}}$. For our model, we consider the chemical reactions (1) and (2). As a detailed simulation of the two-phase flow with resolved interfaces between the phases is very costly, we employ a simpler isothermal and steady-state volume-averaged mixture model which is given by

$$\nabla \cdot (\rho(\alpha)u) = 0 \quad \text{in } \Omega,$$

$$\rho(\alpha)(u \cdot \nabla)u - \nabla \cdot (\mu(\alpha)\nabla u) - \nabla\left(\frac{\mu(\alpha)}{3}\nabla \cdot u\right) + \nabla p + \mu(\alpha)K_0^{-1}u + \rho(\alpha)K_1^{-1}|u|u = 0 \quad \text{in } \Omega, \tag{3}$$

$$\nabla \cdot (\rho^i u\, \alpha^i) - \nabla \cdot (D^i\nabla\alpha^i) = \dot{V}^i \quad \text{in } \Omega \text{ for } i = 1, \dots, N,$$

$$\nabla \cdot \left(\rho^i u\, \alpha^i Y_j^i\right) - \nabla \cdot \left(D_j^i \nabla(\alpha^i Y_j^i)\right) = \dot{V}_j^i \quad \text{in } \Omega \text{ for } j = 1, \dots, M, i = 1, \dots, N.$$

Here, $u$ is the (bulk) fluid velocity, $p$ the pressure, $\rho$ denotes the density, $\mu$ is the dynamic viscosity, and $D$ is a diffusion coefficient. Additionally, $K_0^{-1}$ is a linear viscous resistance coefficient and $K_1^{-1}$ is the quadratic viscous resistance coefficient, which have been obtained with the help of detailed simulations of representative sections of the reactors with Ansys® Fluent, 2022 R2 [38]. The superscript $i$ denotes quantities belonging to phase $i$. The volume fraction of phase $i$ is denoted by $\alpha^i$ and $Y_j^i$ is the mass fraction of (chemical) species $j$ in phase $i$. We consider $N$ phases and $M$ chemical species. In case species $j$ does not occur in phase $i$, we drop the respective equation since $Y_j^i = 0$. Finally, we have two source terms for the volume and species given by $\dot{V}^i$ and $\dot{V}_j^i$, respectively.

The density and viscosity are given by volume averages of the phase-specific quantities, i.e., $\rho(\alpha) = \sum_{i=1}^{N} \alpha^i \rho^i$ and $\mu(\alpha) = \sum_{i=1}^{N} \alpha^i \mu^i$. Moreover, we assume that each phase has the same bulk velocity $u$ so that we only consider a single momentum equation.

For the source terms, we assume to have the following set of chemical reactions in the general form

$$\sum_{i=1}^{N} \sum_{j=1}^{M} \nu_{j,i,k}' \mathcal{M}_j^i \to \sum_{i=1}^{N} \sum_{j=1}^{M} \nu_{j,i,k}'' \mathcal{M}_j^i \quad \text{for } k = 1, \dots, N_r, \tag{4}$$

where $\mathcal{M}_j^i$ is a symbol for species $j$ in phase $i$, $\nu_{j,i,k}'$, $\nu_{j,i,k}''$ are the respective (non-negative) forward and backward stoichiometric coefficients for reaction $k$, and $N_r$ is the number of reactions.

The source terms for the chemical reactions are of the form

$$\dot{V}_j^i = M_j \alpha^i \sum_{k=1}^{N_r} \nu_{j,i,k} Q_k^i, \tag{5}$$

where $\nu_{j,i,k} = \nu_{j,i,k}'' - \nu_{j,i,k}'$ and $M_j$ is the molar weight of species $j$. Due the consistency, the source term for volume fraction $i$ is of the form

$$\dot{V}^i = \sum_j \dot{V}_j^i. \tag{6}$$

The source terms for the volume fractions satisfy

$$\sum_{i=1}^{N} \dot{V}^i = 0, \tag{7}$$

due to the conservation of mass by the chemical reactions. Here, $Q_k^i$ is the rate of reaction $k$ in phase $i$, which is given by

$$Q_k^i = (k_f)_k \prod_{i=1}^{N} [X_j^i]^{\nu_{j,i,k}'}, \tag{8}$$

where $[X_j^i]$ is the molar concentration of species $j$ in phase $i$, $\nu_{j,i,k}'$ is the respective forward stoichiometric coefficient, and $(k_f)_k$ is the $k$-th forward rate coefficient, given by the Arrhenius model

$$(k_f)_k = A_k\, exp\left(-\frac{(E_a)_k}{RT}\right). \tag{9}$$

The concentration $[X_j^i]$ can be computed as

$$[X_j^i] = \frac{\rho^i Y_j^i}{M_j}. \tag{10}$$

We note that this model is sufficient for our case as all educts are part of the liquid phase. In more complex situations, the rate of progress for the reactions may be more complicated.

The above system of partial differential equations (PDEs) is supplemented with the boundary conditions

$$u = u_{\text{in}}, \quad \alpha^i = \alpha_{\text{in}}^i, \quad Y_j^i = (Y_j^i)_{\text{in}} \quad \text{on } \Gamma_{\text{in}},$$

$$u \cdot n = 0, \quad \mu \partial_n u \times n = 0, \quad D^i \partial_n \alpha^i = 0, \quad D_j^i \partial_n(\alpha^i Y_j^i) = 0 \quad \text{on } \Gamma_{\text{wall}}, \tag{11}$$

$$\mu\, \partial_n u - p n = 0, \quad D^i\, \partial_n \alpha^i = 0, \quad D^i_j\, \partial_n\big(\alpha^i Y^i_j\big) = 0 \quad \text{on } \Gamma_{\text{out}}.$$

Hence, at the inlet we use Dirichlet conditions specifying the initial velocity, volume, and mass fractions, on the wall boundary we use a slip boundary condition for the velocity as well as homogeneous Neumann conditions for volume and mass fractions, and on the outlet, we use a do-nothing condition for velocity and pressure as well as homogeneous Neumann conditions for volume and mass fractions.

In our setting, we have $N = 2$ phases as we consider a liquid as well as a gas phase and use the superscripts "gas" and "liquid" to distinguish between them. Moreover, we consider $M = 7$ chemical species.

Note that we use an effective one-dimensional model of the reactor for our numerical experiments, substantially reducing the numerical effort for solving our model. This is achieved by averaging the model presented above over the cross section of the reactors. We refer the reader to [39] where a similar approach is used. For the numerical solution of the model we use the finite element software FEniCS [40].

### 2.5 Grey-Box PFR Model

Our second model involves a one-dimensional PFR-like [41] description of the continuous DAN formation. It is described by the following initial value problem:

$$
\begin{aligned}
& \frac{d[X_j]}{dz} = \frac{1}{u}\, \alpha^{\text{liquid}} \sum_{k=1}^{N_r} \nu_{j,k} Q_k, \\
& \rho u = \big(\rho^{\text{liquid}}\alpha^{\text{liquid}} + \rho^{\text{gas}}\alpha^{\text{gas}}\big) u = \text{const}, \\
& \alpha^{gas} = \frac{[X_{N_2}]M_{N_2}}{\rho^{\text{gas}}} \quad \text{and} \quad \alpha^{\text{liquid}} = 1 - \alpha^{\text{gas}}, \\
& [X_j] = [X_j]_{\text{in}} \quad \text{on } \Gamma_{\text{in}}.
\end{aligned}
\tag{12}
$$

In contrast to the physical one, the grey-box model uses the concentrations $[X_j]$ of species $j$ w.r.t. the entire volume, and not the phase-specific one. The notation is same as in the physical model. We assume the gas density is given by the ideal gas law

$$\rho^{gas} = \frac{p M_{N_2}}{RT}, \tag{13}$$

where we suppose the pressure to be uniform along the reactor. Note that the PFR model 12 is significantly less complex compared to the physical model (3), as it, e.g., neglects the coupling of pressure and velocity, no momentum equation is solved, and diffusion effects are neglected. For these reasons, it is cheaper to be solved numerically and, thus, better suited for introducing a data-based component [42].

Whereas the source term for the DAN synthesis $Q_1$ takes the previously considered form given in Eq. (8), a data-driven function of the form

$$Q_2 = f([X_{DAN}], T, p)\, [X_{DAN}] \tag{14}$$

is used as source term for the DAN decomposition (2). Here, $f$ denotes a neural network consisting of a single hidden layer with 4 neurons. The activation function is chosen to be $\max(0, x^3)$ and it is used after the hidden and output layers, guaranteeing non-negative reactions rates [43].

### 2.6 Identification of Parameters for Reaction Kinetics

The reaction kinetics for the DAN synthesis and decomposition are not known in the literature. For this reason, we now aim at numerically identifying these parameters by calibrating our models to the experimental data.

### 2.6.1 Cost Function for the Mixer

In case of the micro-structured reactor, FTIR measurements, proportional to the molar concentration of DAN, are taken. However, the proportionality factor $\gamma$, relating the measured band areas to the concentration of DAN, is unknown and, thus, has to be fitted numerically. For these reasons, we introduce the following least-squares cost function, to be minimized

$$J_1 = \tfrac{1}{2} \sum_{m=1}^{N_{\text{mix}}} \int_{\Gamma_{\text{out}}} \left( B_m^{\text{sim}} - B_m^{\text{exp}} \right)^2 \mathrm{d}s, \tag{15}$$

where $B_m^{\text{sim}}$ and $B_m^{\text{exp}}$ denote the simulated and measured DAN band areas for experiment $m$, respectively, where the former is given by

$$B_m^{\text{sim}} (y, u) = \frac{\left[ x_{\text{DAN}}^{\text{liquid}} \right]}{\gamma}, \tag{16}$$

and the concentration of DAN can be computed with the help of Eq. (10).

### 2.6.2 Cost Functions for the Calorimeter

#### 2.6.2.1 Heat Released by the Reactions

In case of the calorimeter, we first focus on the spatially resolved heat flux measurements. Let $\Delta H_1$ and $\Delta H_2$ denote the enthalpies of reactions (1) and (2), respectively. Moreover, let $\Omega_i$ be the $i$-th zone of the reactor geometry $\Omega$. Then, the heat released by the reaction in zone $i$ is given by

$$q_i^{\text{sim}} = \int_{\Omega_i} \alpha^{\text{liquid}} \Delta H_1 Q_1 + \alpha^{\text{liquid}} \Delta H_2 Q_2 \; \mathrm{d}x, \tag{17}$$

where $Q_i = Q_i^{\text{liquid}} + Q_i^{\text{gas}}$. We introduce the following least-squares cost function

$$J_2 = \tfrac{1}{2} \sum_{m=1}^{N_{\text{calo}}} \sum_{i=1}^{N_{\text{zones}}} \left( q_i^{m,\text{sim}} - q_i^{m,\text{exp}} \right)^2, \tag{18}$$

where $q_i^{m,\text{sim}}$ and $q_i^{m,\text{exp}}$ denote the simulated and measured values of the heat released by the reaction (cf. eq. (17)) and $N_{\text{calo}}$ and $N_{\text{zones}} = 7$ are the number of experiments for the calorimeter and the number of zones, respectively.

#### 2.6.2.2 Gas Flow Rate

Second, we consider the gas flow rate measurements. The gas flow rate can be computed via

$$G^{\text{sim}} = \frac{1}{\rho_{\text{ref}}^{\text{gas}}} \int_{\Gamma_{\text{out}}} \alpha^{\text{gas}} \rho^{\text{gas}} u \cdot n \; \mathrm{d}s, \tag{19}$$

where the factor $\frac{1}{\rho_{\text{ref}}^{\text{gas}}}$ is used to convert the simulated flow rate to the reference thermodynamic conditions, i.e., $\rho_{\text{ref}}^{\text{gas}}$ is the density of $N_2$ at $p = 1$ atm and $T = 20\ ^{\circ}\text{C}$. Hence, we introduce the following cost function

$$J_3 = \tfrac{1}{2} \sum_{m=1}^{N_{\text{calo}}} \left( G_m^{\text{sim}} - G_m^{\text{exp}} \right)^2, \tag{20}$$

that measures the discrepancy between simulated ($G_m^{\text{sim}}$) and measured ($G_m^{\text{exp}}$) gas flow rates.

### 2.6.3 Numerical Solution of Parameter Identification Problems

In the previous sections, we have introduced three cost functions measuring the least-squares error between our models and experimental measurements for the two reactors. To minimize them simultaneously, we use a weighted-sum scalarization (each cost function is additionally normalized by a scalar $J_{k,0}$ to make them comparable in magnitude) [44]

$$J = \lambda_1 \frac{J_1}{J_{1,0}} + \lambda_2 \frac{J_2}{J_{2,0}} + \lambda_3 \frac{J_3}{J_{3,0}} \tag{21}$$

which leads to the following parameter identification problems

$$\min_{y_{\text{phys}}, u_{\text{phys}}} J \quad \text{subject to the physical model (3),} \tag{22}$$

and

$$\min_{y_{\text{PFR}}, u_{\text{PFR}}} J \quad \text{subject to the grey-box PFR model (12).} \tag{23}$$

Here, $y_{\text{phys}} = [u, p, \alpha, Y]^T$ denotes the state variables of the physical model, where $\alpha$ is a vector containing the volume fractions $\alpha_i$ and $Y$ is a vector containing the species mass fractions $Y_j^i$. The control variables for this model are given by $u_{\text{phys}} = [\log(A_1), (E_a)_1, \log(A_2), (E_a)_2, \gamma, \Delta H_1, \Delta H_2]^T$. Moreover, for the grey-box PFR model, $y_{\text{PFR}} = X$ is the vector of species concentrations $[X_j^i]$ and the control variables are given by $u_{\text{PFR}} = [\log(A_1), (E_a)_1, \gamma, \Delta H_1, \Delta H_2, W, b]^T$, where $W$ and $b$ denote the weights and biases, respectively, of the neural network appearing in the data-driven component (cf. Section 2.5). Note that the grey-box PFR model introduces considerably more model parameters than the physical one, as 16 weights and 5 biases are present in the neural network. Moreover, we use the (natural) logarithm of the pre-exponential factors as control variables as this ensures a better scaling (cf. [39]).

We choose the following weights $\lambda$ for the cost functions: As the physical model is currently unable to model the gas flow rate accurately, we only consider cost functions $J_1$ and $J_2$, so we adopt a weight vector of $\lambda = [1, 1, 0]^T$. The grey-box model, on the other hand, has been fitted considering all data with a corresponding weight vector $\lambda = [1, 1, 1]^T$.

For the physical model, we solve the parameter identification problem (22) numerically with the software package cashocs [45, 46], which automates the solution of such PDE constrained problems using a discretization of a continuous adjoint approach and is developed at Fraunhofer ITWM. For the sake of brevity, we refer the reader to [45, 46], where cashocs is presented, and to [39], where it is used to solve a similar parameter identification problem in the context of a gas phase reaction. The grey-box model is implemented in the dynamic programming language Julia [47] and IPOPT [48] is used to solve the parameter identification problem (23).

Due to the high number of model parameters in the grey-box PFR model, we use $k$-fold cross-validation to check overfitting is not occurring [49]. This works as follows: given a total of $k$ data sets (in our case experimental data for the mixer and calorimeter setups), we exclude one data set from the identification and solve (23) with the remaining $k - 1$ data sets. The excluded data set is then used as test set, i.e., we evaluate the model prediction and compare it with the actual data. This is then repeated $k$ times so that each case is used as test case. With a total of 24 and 33 data sets for mixer and calorimeter, respectively, this means (23) is solved 57 times. To exclude that overfitting occurred, the average model prediction error for the test is supposed to be comparable to the error achieved for the training set. We use the mean average error (MAE) to quantify the prediction quality of the model [50].

## 3 Results and Discussion

Now, we discuss the experimental results and analyze the quality of the calibrated simulations for both models.

### 3.1 Mixer – FTIR Measurements

We start with analyzing the DAN formation by means of the band areas of the characteristic DAN absorptions from the mixer setup (Figure 4). In general, for equal residence times, higher temperatures cause growing values of the band area, i.e., DAN concentration. Moreover, at low temperatures (20 and 35 °C), we observe a monotonic behavior of the band area w.r.t. the residence time. However, the experiments performed at 55 and 70 °C show a non-monotonic dependence, with a clear maximum achieved at residence time values between 60 and 80 s. Apparently, at higher temperatures, both the DAN generation and decomposition occur more rapidly. The gas content in the system increases and consequently the flow velocity in the reactor changes significantly. This leads to fluctuations and an overall reduction in the residence time that directly affects the DAN generation. Thus, the role of decomposition and a high presence of gas are supposed to be the main reasons for the decrease of DAN at large residence times.

There is a good agreement between predictions of both models and data. In particular, the physical model captures the experimental trends well both quantitatively and qualitatively, and even the non-monotonic behavior w.r.t. the residence time can be seen. In contrast, the grey-box model is less detailed, particularly at higher temperatures and longer residence times. Moreover, the grey-box model is not able to capture the non-monotonic behavior described earlier, due to the significant amount of gas, which strongly limits both reactions. This can be also seen in the coefficient of determination, which is lower for the grey-box model at $R^2 = 0.88$ compared to $R^2 = 0.95$ for the physical model (Figure 5).

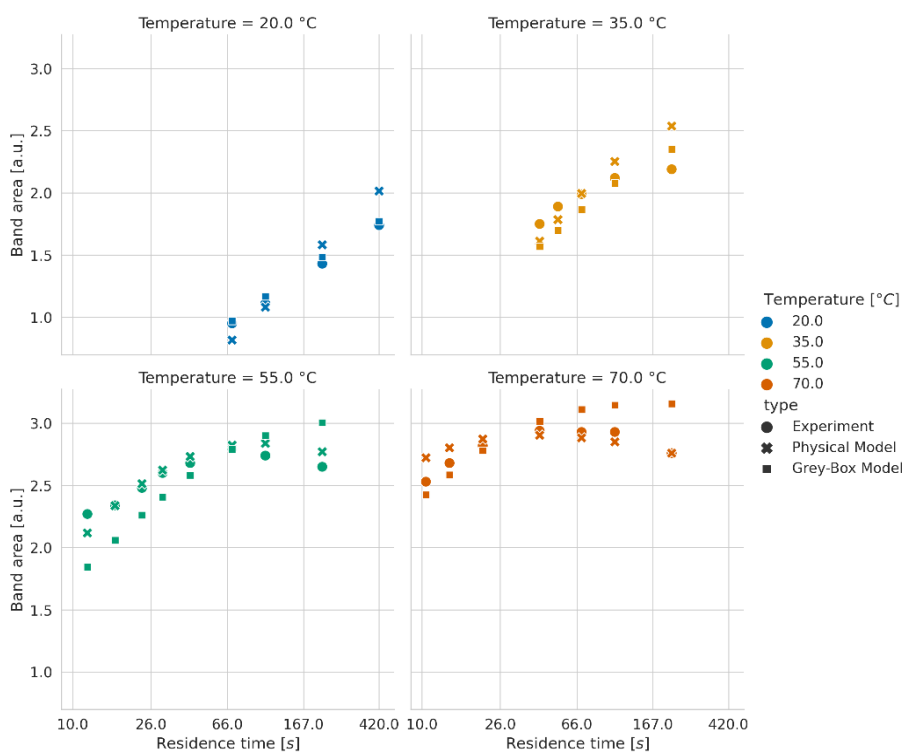

**Figure 4.** Comparison of simulated and measured FTIR absorption band areas over the residence time in the mixer for various temperatures.

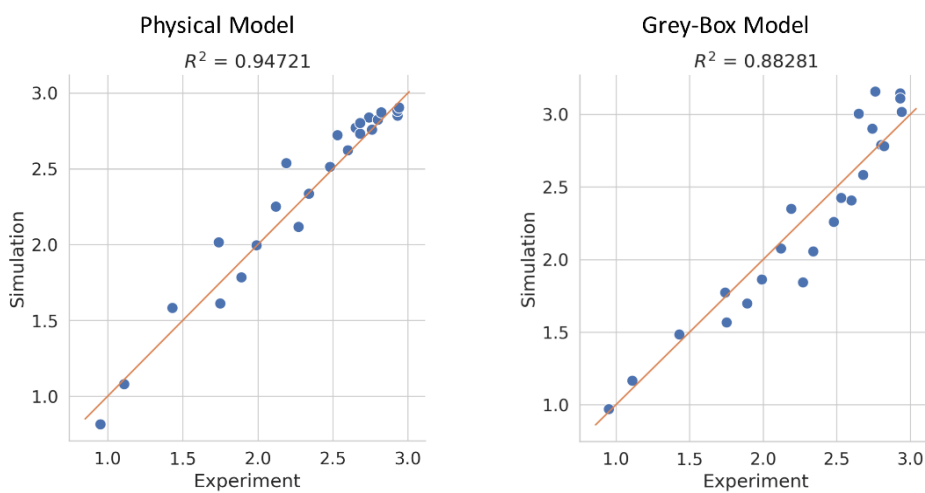

**Figure 5.** Parity plots comparing the simulated and measured FTIR absorption band areas for the mixer (left - physical model, right - grey-box PFR model).

### 3.2 Calorimeter

#### 3.2.1 Heat Released by the Reactions

From the experimental data, we observe that the total (i.e., sum of the spatially resolved values) measured heat released by reactions (1) and (2) increases with temperature and flow rate (Figure 6). First, as expected, a higher temperature results in a higher reaction rate. Second, a higher flow rate implies that more educts are fed into the system and that they have less time to react, resulting in a higher reactivity. Finally, we observe that there is little dependence of the released heat on the pressure as there is only a minor increase with growing pressures.

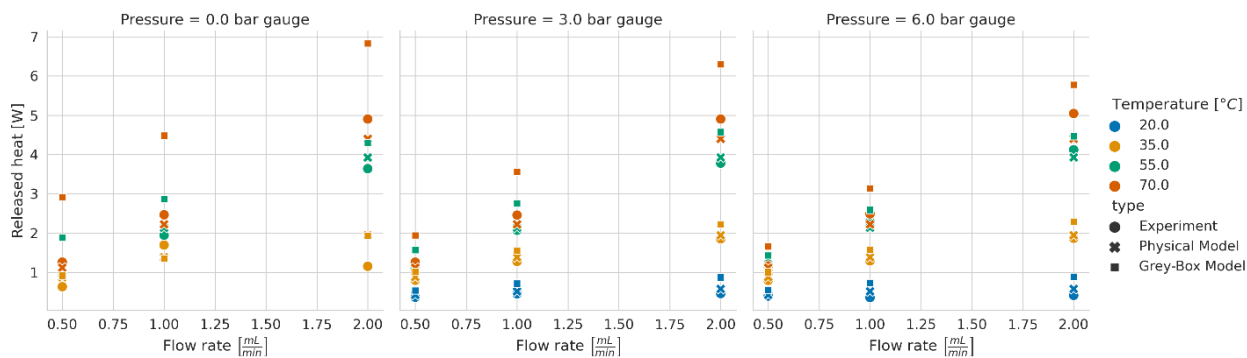

**Figure 6.** Comparison between simulation and experiments of the total heat released by the reaction over the flow rate in the reaction calorimeter for various gauge pressures and temperatures.

Investigating the heat released at each zone, which is shown in Figure 7, we observe that it decays approximately exponentially along the reaction channel. Particularly, most of the heat is released in the first zone, which suggests that the reactions occur very rapidly in the beginning of the calorimeter's mixing channel and are then slow down significantly throughout the remaining zones.

As can be seen from Figure 6 and Figure 7, both models predict the measured data well. In particular, the coefficient of determination satisfies $R^2 = 0.98$ for the physical and $R^2 = 0.92$ for the grey-box model, indicating a very good fit (Figure 8). As for the FTIR data discussed previously, the physical model achieves a better fit than the grey-box model. We note that the physical model achieves an excellent fit for the released heat, even though it does only consider a rather slow decomposition reaction (cf. Section 3.3). This indicates that there may be multiple parameter configurations for which similar heat released by the combination of DAN synthesis and decomposition (and other parallel reactions not considered) can be achieved with our simulation models.

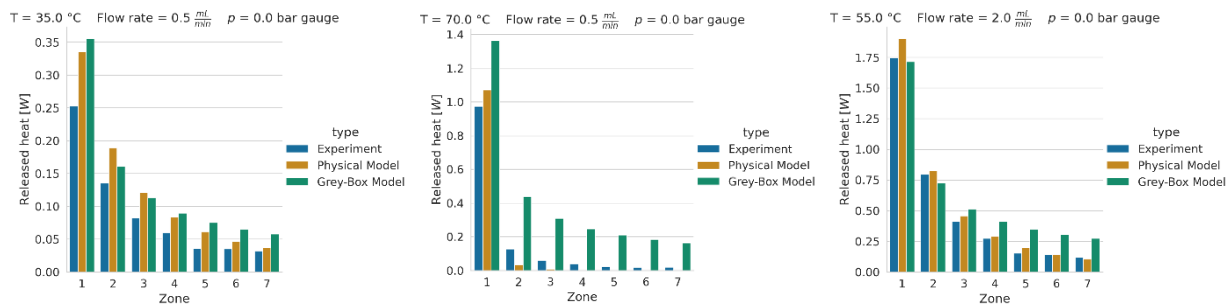

**Figure 7.** Comparison of measured and simulated heat released by the reaction over the measurement zones of the flow calorimeter, shown for three exemplary cases. Here, $\tau$ denotes the residence time.

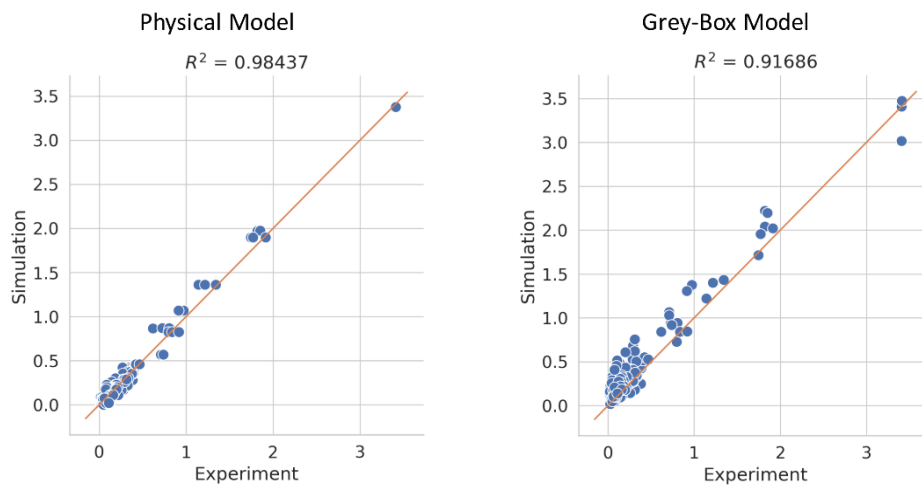

**Figure 8.** Parity plots comparing the simulated and measured heat released by the reactions in each calorimeter zone (left - physical model, right - grey-box PFR model).

### 3.2.2 Gas Flow Rate

Finally, we discuss the gas flow rates, which are shown in Figure 9. We observe that, as expected, the gas flow rate increases with temperature and that it does not depend strongly on the pressure, similarly to the released heat. Moreover, we observe that the gas flow rate increases almost linearly with the flow rate. This is expected as lower flow rates imply lower amounts of educts, so that less overall gas can be formed.

As described in Section 2.4, we could not calibrate the physical model to the gas flow rate data, yet. With the identified kinetic parameters from Table 1, the physical model strongly under-predicts the measured values. Further model refinements are object of future work.

These results originally motivated our choice of employing a grey-box model, where a data-driven component is used to model the DAN decomposition. From the results shown in Figure 9, we observe that the grey-box model can  qdescribe the gas flow rate accurately. Moreover, the parity plot shown in Figure 10 shows a very good agreement between measurements and grey-box model, with a corresponding coefficient of determination of $R^2 = 0.96$.

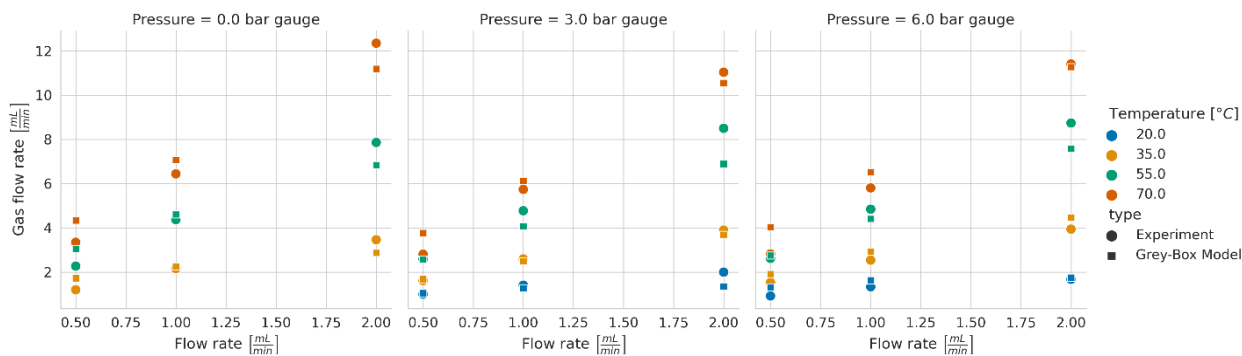

**Figure 9.** Comparison between experiments and grey-box PFR predictions of the gas flow rate (evaluated at atmospheric pressure and 20 °C) over the flow rate for various temperatures and gauge pressures.

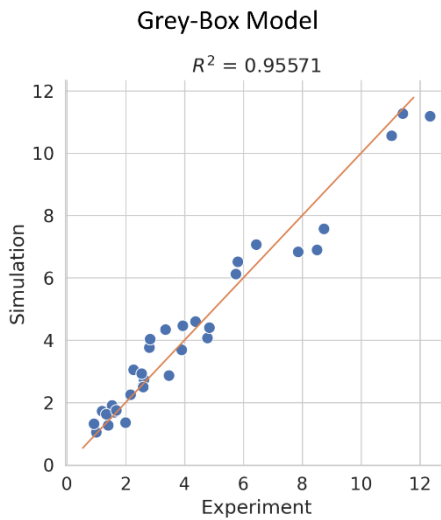

**Figure 10.** Parity plot comparing the simulated and measured gas flow rates per educt for the reaction calorimeter for the grey-box PFR model.

### 3.3 Obtained Reaction Parameters

The reaction parameters obtained with our parameter identification approach are shown in Table 1. Note, that the parameters for $(E_a)_2$ and $\log(A_2)$ are not shown for the grey-box model as we employed a neural network to predict the source term for this model (cf. Section 2.5). We observe that both the physical and grey-box model yield similar parameters for the DAN synthesis reaction (1) with an activation energy of $80 \frac{\text{kJ}}{\text{mol}}$ and $73 \frac{\text{kJ}}{\text{mol}}$ and a value for $\log(A)$ of 21 and 18, respectively. Moreover, even the values for the proportionality factor of 339 and 295 as well as the obtained values for the enthalpy of the DAN synthesis of $123 \frac{\text{kJ}}{\text{mol}}$ and $136 \frac{\text{kJ}}{\text{mol}}$ are in good agreement between the models. This indicates that both models predict the synthesis of DAN in a rather similar fashion.

Due to the fact that the physical model does neglect the measured gas flow rates (cf. Section 2.6.3) and the grey-box model does not, there are some (expected) differences between the models w.r.t. the DAN decomposition reaction (2). In particular, the physical model predicts a vanishing activation energy, so that the reaction would not be temperature dependent, and a very low pre-exponential factor, so that the DAN decomposition is rather weak. In fact, we are not able to fit the gas flow rates with the physical model. For this reason, the reaction parameters corresponding to the DAN decomposition reaction (2), namely $(E_a)_2$, $\log(A_2)$, and $\Delta H_2$, are not reliable and are most likely not physically reasonable.

In contrast, for the grey-box model, the source term for the DAN decomposition, which is shown in Figure 11, shows a clear temperature dependence. As expected, higher temperatures yield higher reaction rates. Additionally, the source term is also pressure-dependent, where higher values are achieved at high pressures. Finally, the source depends non-linearly and even non-monotonically on the DAN concentration. This could be reasonable as the DAN decomposition reaction (2) and possible subsequent reactions are not modeled explicitly. However, the influence of such reactions may be present in the measured data that are accessible to our grey-box model.

**Table 1.** Reaction kinetics parameters identified with the physical and grey-box PFR model.

| Parameter | Name | physical | Grey-box PFR | Unit |
|-----------|------|----------|--------------|------|
| $(E_a)_1$ | Activation energy for DAN synthesis | 80.25 | 72.57 | $\frac{kJ}{mol}$ |
| $log(A_1)$ | Logarithm of pre-exponential factor for DAN synthesis | 20.6 | 18.1 | — |
| $(E_a)_2$ | Activation energy for DAN decomposition | 0.0[a)] | -[b)] | $\frac{kJ}{mol}$ |
| $log(A_2)$ | Logarithm of pre-exponential factor for DAN decomposition | -10.5[a)] | -[b)] | — |
| $\gamma$ | Proportionality factor for the band area | 339.48 | 295.45 | Dimensionless |
| $\Delta H_1$ | Enthalpy of DAN synthesis | 122.7 | 135.7 | $\frac{kJ}{mol}$ |
| $\Delta H_2$ | Enthalpy of DAN decomposition | 102.2[a)] | 26.0 | $\frac{kJ}{mol}$ |

a) obtained parameter is not reliable; b) obtained with neural network.

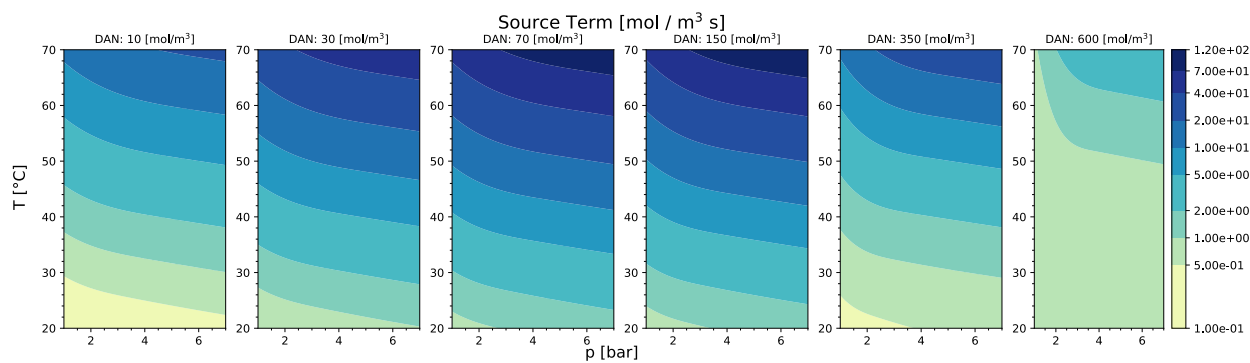

**Figure 11.** Data-driven source term (DAN decomposition) for the grey-box PFR model evaluated at different temperatures, pressures, and DAN concentrations.

Results concerning $k$-fold cross-validation for the grey-box model are shown in Table 2. The results do not emphasize overfitting as the average MAE for test and training set is comparable for all measured quantities.

**Table 2.** $k$-fold cross-validation for the grey-box PFR model. Average MAE on training and test sets.

| Set | DAN [a.u.] | Heat [W] | Gas Flow Rate [ml/min] |
|-----|-----------|----------|------------------------|
| Training | 0.193 | 0.404 | 0.0837 |
| Test | 0.259 | 0.516 | 0.0903 |

## 4  Conclusion and Outlook

In this paper, the synthesis of diazo acetonitrile (DAN) in continuous flow reactors has been experimentally investigated at two setups. In the first, a micro-structured reactor with a chaotic mixing structure, DAN band areas are measured by FTIR spectroscopy. The second reactor, a flow calorimeter, provides spatially resolved information on the heat released by the reactions. Additionally, the gas flow rate of nitrogen (formed as a byproduct in the DAN decomposition) is measured.

We presented two numerical models for the investigated continuous flow reactors: First, a two-phase physical model of the chemically reactive flow, and second, a grey-box plug flow reactor model using techniques from machine learning. We show that both models can be calibrated well with experimental data for DAN formation. Due to the high complexity of the reactions involved in the synthesis (and decomposition) of DAN, it is very complex to derive a suitable model of the chemically reactive flow that incorporates all important physical phenomena. However, we have seen that grey-box modeling can be used to greatly simplify the modeling efforts, while also obtaining comparable, and for the case of the gas flow rates, even superior results. A combination of the advantages of both approaches is up for future work.

The results obtained in this paper give detailed insights into the chemical and physical processes occurring during the production of DAN and are building blocks for a further investigation of the DAN synthesis. Particularly, the further development of the DAN synthesis with an additional extraction of the DAN to an organic phase, where it is more stable, is of major importance and currently under investigation, as is, of course, adapting our simulation models to this extended setting.

**Supporting Information**

Supporting Information for this article can be found under [Link provided by Wiley]. This section includes additional references to primary literature relevant for this research [XX–XY].

**Acknowledgment**

This work was funded by the Fraunhofer Lighthouse project ShaPID.

**Symbols used**

$A_1$ $\quad\left[\frac{1}{s}\left(\frac{mol}{m^3}\right)^{-1}\right]$ $\quad$ pre-exponential factor for DAN synthesis

$A_2$ $\quad$ [1/s] $\quad$ pre-exponential factor for DAN decomposition

$B$ $\quad$ [a.u.] $\quad$ band area of DAN

$b$ $\quad$ [-] $\quad$ biases for the neural network

$D$ $\quad\left[\frac{kg}{ms}\right]$ $\quad$ diffusion coefficient

$(E_a)_1$ $\quad$ [J/mol] $\quad$ activation energy for the DAN synthesis

$(E_a)_2$ $\quad$ [J/mol] $\quad$ activation energy for the DAN decomposition

$f$ $\quad$ [1/s] $\quad$ data-driven function for modeling the DAN decomposition in the grey-box model

$G$ $\quad$ [mL/min] $\quad$ Gas flow rate of nitrogen

$\Delta H_1$ $\quad$ [J/mol] $\quad$ enthalpy of the DAN synthesis

$\Delta H_2$ $\quad$ [J/mol] $\quad$ enthalpy of the DAN decomposition

$J$ $\quad$ [-] $\quad$ cost function for the parameter identification

$K_0$ $\quad$ [m²] $\quad$ viscous permeability

$K_1$ $\quad$ [m] $\quad$ inertial permeability

$(k_f)_1$ $\quad\left[\frac{1}{s}\left(\frac{mol}{m^3}\right)^{-1}\right]$ $\quad$ forward rate coefficient for DAN synthesis

$(k_f)_2$ $\quad$ [1/s] $\quad$ forward rate coefficient for DAN decomposition

$\log(x)$ $\quad$ [-] $\quad$ natural logarithm of $x$

$M$ $\quad$ [-] $\quad$ number of species

| $M_j$ | [kg/mol] | molar weight of species $j$ |
|---|---|---|
| $\mathcal{M}_j^i$ | [-] | symbol for species $j$ in phase $i$ |
| $N$ | [-] | number of phases |
| $N_{\text{calo}}$ | [-] | number of experiments considered for the calorimeter setup |
| $N_{\text{mix}}$ | [-] | number of experiments considered for the mixer setup |
| $N_{\text{zones}}$ | [-] | number of distinct zones for the reaction calorimeter |
| $p$ | [Pa] | pressure |
| $q$ | [W] | heat released by the reactions |
| $Q_k^i$ | [mol/m³/s] | reaction rate in phase $i$ for reaction $k$ |
| $R$ | $\left[\frac{\text{J}}{\text{K mol}}\right]$ | universal gas constant |
| $R^2$ | [-] | coefficient of determination |
| $T$ | [K] | temperature |
| $u$ | [m/s] | flow velocity |
| $u_{\text{PFR}}$ | [-] | control variables for the grey-box PFR model |
| $u_{\text{phys}}$ | [-] | control variables for the physical model |
| $\dot{V}^i$ | $\left[\frac{\text{kg}}{\text{m}^3\,\text{s}}\right]$ | source term for volume of phase $i$ |
| $\dot{V}_j^i$ | $\left[\frac{\text{kg}}{\text{m}^3\,\text{s}}\right]$ | source term for species $j$ of phase $i$ |
| $W$ | [-] | weights for the neural network |
| $X$ | [mol/m³] | vector containing the molar concentrations $[X_j]$ |
| $[X_j^i]$ | [mol/m³] | molar concentration of species $j$ in phase $i$ |
| $Y_j^i$ | [-] | mass fraction of species $j$ in phase $i$ |
| $y_{\text{phys}}$ | [-] | state variables for the physical model |
| $y_{\text{PFR}}$ | [-] | state variables for the grey-box PFR model |

**Greek letters**

| $\alpha$ | [-] | volume fraction |
|---|---|---|
| $\Gamma$ | [-] | boundary of $\Omega$ |
| $\gamma$ | [-] | proportionality factor between DAN band area and concentration |
| $\lambda$ | [-] | weight vector for scalarization of the cost functions |
| $\mu$ | [Pa s] | dynamic viscosity |
| $\nu_{j,i,k}$ | [-] | difference of backward and forward stoichiometric coefficients |
| $\nu'_{j,i,k}$ | [-] | forward stoichiometric coefficient of species $j$ in phase $i$ for reaction $k$ |
| $\nu''_{j,i,k}$ | [-] | backward stoichiometric coefficient of species $j$ in phase $i$ for reaction $k$ |
| $\rho$ | [kg/m³] | density |
| $\Omega$ | [-] | computational domain of the chemically reactive flow |

| $\Omega_i$ | [-] | $i$-th zone of the reaction calorimeter |

**Sub- and Superscripts**

| exp | experimental |
| gas | belonging to the gaseous phase |
| $i$ | belonging to phase $i$ |
| in | belonging to the inlet |
| $j$ | belonging to species $j$ |
| $k$ | belonging to reaction $k$ |
| liquid | belonging to the liquid phase |
| out | belonging to the outlet |
| ref | reference |
| sim | simulated |
| wall | belonging to the wall boundary |

**Abbreviations**

ATR – attenuated total reflection

DAN – diazo acetonitrile

FTIR spectrometer – Fourier-transform infrared spectrometer

MAE – mean average error

PFR – Plug-Flow Reactor

Supporting Information

# Continuous Synthesis of Diazo Acetonitrile: Experiments, Modeling, and Parameter Identification

Dr. Marco Baldan[1,*], Dr. Sebastian Blauth[1], Dr. Dušan Bošković[2], Dr. Christian Leithäuser[1], Dr. Alexander Mendl[2], Ligia Radulescu[2], Maud Schwarzer[2], Heinrich Wegner[2], and Prof. Dr. Michael Bortz[1]

*Correspondence:* marco.baldan@itwm.fraunhofer.de, Fraunhofer Institute for Industrial Mathematics ITWM, Fraunhofer-Platz 1, 67663 Kaiserslautern, Germany

## 5    Additional Results and Discussions

In the following, we provide the results obtained by our experiments for the mixer and calorimeter setups.

**Table S1.** Overview of the results from the online FT-IR measurement

| Residenz time [s] | Temperature [°C] | Band area [a.u.] |
|---|---|---|
| 10,5 | 70 | 2,53 |
| 12,0 | 55 | 2,27 |
| 14,0 | 70 | 2,68 |
| 16,8 | 55 | 2,34 |
| 21,0 | 70 | 2,82 |
| 23,3 | 55 | 2,48 |
| 30,0 | 55 | 2,6 |
| 42,0 | 35 | 1,75 |
| 42,0 | 55 | 2,68 |
| 42,0 | 70 | 2,94 |
| 52,5 | 35 | 1,89 |
| 70,0 | 19 | 0,95 |
| 70,0 | 35 | 1,99 |
| 70,0 | 55 | 2,8 |
| 70,0 | 70 | 2,93 |
| 105,0 | 19 | 1,11 |
| 105,0 | 35 | 2,12 |
| 105,0 | 55 | 2,74 |
| 105,0 | 70 | 2,93 |
| 210,0 | 19 | 1,43 |
| 210,0 | 35 | 2,19 |
| 210,0 | 55 | 2,65 |

| | | | | |
|---|---|---|---|---|
| 210,0 | 70 | 2,76 | | |
| 420,0 | 20 | 1,74 | | |

**Table S2.** Overview of the results from the heat flow measurement

| Residence time [s] | Temperature [°C] | Pressure gauge [bar] | Heat [W] | Gas flow rate [mL/min] |
|---|---|---|---|---|
| 26,7 | 35 | 0 | 0,63 | 1,63 |
| 26,7 | 55 | 0 | 0,94 | 2,72 |
| 26,7 | 70 | 0 | 1,26 | 4,44 |
| 26,7 | 20 | 3 | 0,37 | 1,00 |
| 26,7 | 35 | 3 | 0,80 | 1,41 |
| 26,7 | 55 | 3 | 1,09 | 2,00 |
| 26,7 | 70 | 3 | 1,26 | 0,92 |
| 26,7 | 20 | 6 | 0,40 | 1,34 |
| 26,7 | 35 | 6 | 0,80 | 1,67 |
| 26,7 | 55 | 6 | 1,21 | 1,20 |
| 26,7 | 70 | 6 | 1,16 | 2,16 |
| 53,3 | 35 | 0 | 1,69 | 3,46 |
| 53,3 | 55 | 0 | 1,95 | 1,62 |
| 53,3 | 70 | 0 | 2,47 | 2,60 |
| 53,3 | 20 | 3 | 0,45 | 3,90 |
| 53,3 | 35 | 3 | 1,28 | 1,53 |
| 53,3 | 55 | 3 | 2,07 | 2,54 |
| 53,3 | 70 | 3 | 2,46 | 3,94 |
| 53,3 | 20 | 6 | 0,36 | 2,27 |
| 53,3 | 35 | 6 | 1,29 | 4,37 |
| 53,3 | 55 | 6 | 2,29 | 7,86 |
| 53,3 | 70 | 6 | 2,49 | 2,59 |
| 106,7 | 35 | 0 | 1,15 | 4,77 |
| 106,7 | 55 | 0 | 3,64 | 8,50 |
| 106,7 | 70 | 0 | 4,90 | 2,62 |
| 106,7 | 20 | 3 | 0,46 | 4,84 |
| 106,7 | 35 | 3 | 1,85 | 8,74 |
| 106,7 | 55 | 3 | 3,78 | 3,35 |
| 106,7 | 70 | 3 | 4,90 | 6,44 |
| 106,7 | 20 | 6 | 0,41 | 12,35 |
| 106,7 | 35 | 6 | 1,87 | 2,81 |
| 106,7 | 55 | 6 | 4,12 | 5,74 |
| 106,7 | 70 | 6 | 5,04 | 11,04 |